\input amstex
\magnification 1200

\documentstyle{amsppt}
\NoBlackBoxes

\pageheight{22 true cm}
\pagewidth{17.5 true cm}
\vcorrection {-0.0 true cm}
\hcorrection {0 true cm}
\overfullrule=0pt

\pageno=1

\define\R{\text{\rm Re }}

\topmatter \rightheadtext{Optimal control...}
\author
Dmitri Prokhorov \& Alexander Vasil'ev
\endauthor
\title{Optimal control in Bombieri's and \\ Tammi's conjectures}
\endtitle
\address{Department of Mathematics and Mechanics, Saratov State University, Saratov 410026,
Russia}\endaddress
\email{ProkhorovDV\@info.sgu.ru}\endemail
\address{Departamento de Matem\'atica, UTFSM, Casilla, 110-V Valpara\'\i so, Chile}
\endaddress
\email{alexander.vasiliev\@usm.cl}\endemail
\thanks{The first author is supported by the Russian Foundation
for Basic Research, grant \# 01-01-00123 and the grant of the
Ministry of Higher Education (Russia) \# E02-1.0-178; the second
author is supported by FONDECYT (Chile), projects \# 1030373,
\#1040333}, and UTFSM \#12.03.23 \endthanks 
\endkeywords
\abstract Let $S$ stand for the usual class of univalent regular
functions in the unit disk $U=\{z:\,|z|<1\}$ normalized by
$f(z)=z+a_2z^2+\dots$ in $U$, and let $S^M$ be its subclass defined
by restricting $|f(z)|<M$ in $U$, $M\geq 1$. We consider two
classical problems: Bombieri's coefficient problem for the class $S$
and the sharp estimate of the fourth coefficient of a function from
$S^M$. Using L\"owner's parametric representation and the optimal
control method we give exact initial Bombieri's numbers and derive a
sharp constant $M_0$, such that for all $M\geq M_0$ the Pick
function gives the local maximum to $|a_4|$. Numerical approximation
is given.
\endabstract
\endtopmatter
\document

\head{\S 1. Introduction}\endhead

Let $S$ stand for the class of all holomorphic and univalent
functions $f(z)=z+\sum_{n=2}^{\infty}a_nz^n$ in the unit disk
$U=\{z:\,|z|<1\}$. Its subclass of bounded maps $|f(z)|<M$, $M\geq
1$, we denote by $S^M$, $S^{\infty}\equiv S$. During the long
history of univalent functions the famous Koebe function $$
K(z)=\frac{z}{(1-z)^2}=\sum_{n=1}^{\infty}nz^n\in S $$ has been
known to be extremal in many problems. A relevant sample is the
most celebrated Bieberbach Conjecture [2] $|a_n|\leq n$ proved by
L. de Branges in 1984 [5,6]. In spite of many works about
coefficient estimates in the class $S$, there are some difficult
problems that are still unsolved, in particular, the Bombieri
problem and the sharp upper bound for $|a_n|$, $n\geq 4$, for the
subclass $S^M$ that we will deal with.

E.~Bombieri [4] in 1967 posed the problem to find
$$
\sigma_{mn}:=\liminf\Sb f\to K\\f\in S\endSb\frac{n-\R a_n}{m-\R a_m},
\;\;\;m,n\geq2,
$$
$f\to K$ locally uniformly in $U$. We call $\sigma_{mn}$ the Bombieri numbers.
 He conjectured that $\sigma_{mn}=B_{mn}$, where
$$
B_{mn}=\min\limits_{\theta\in [0,2\pi)}\frac{n\,\sin \theta-\sin(n\theta)}{m\,\sin \theta-\sin(m\theta)}.
$$
and proved that $\sigma_{mn}\leq B_{mn}$ for $m=3$ and $n$ odd. It
is noteworthy that D.~Bshouty and W.~Hengartner [7] proved
Bombieri's conjecture for functions from $S$ having real
coefficients in their Taylor's expansion. Continuing this
contribution by D.~Bshouty and W.~Hengartner, the conjecture for
the whole class $S$ has been recently disproved by R.~Greiner and
O.~Roth [9] for $n=2$, $m=3$, $f\in S$. Actually, they have got
the sharp Bombieri number $\sigma_{32}=(e-1)/4e<1/4=B_{32}$.

It is easily seen that $\sigma_{43}=B_{43}=\sigma_{23}=B_{23}=0$.
Applying L\"owner's parametric representation for univalent
functions and the optimal control method we will find the exact
Bombieri numbers $\sigma_{42},\sigma_{24}, \sigma_{34}$ and their
numerical approximations $\sigma_{42}\approx 0.050057\dots$,
$\sigma_{24}\approx 0.969556\dots$, and $\sigma_{34}\approx
0.791557\dots$ (the Bombieri conjecture for these permutations of
$m,n$ suggests $B_{42}=0.1$, $B_{24}=1$, $B_{34}=0.828427\dots$). Of
course, our method permits us to reprove the result of [9] about
$\sigma_{32}$.

Our next target is the fourth coefficient $a_4$ of a function from
$S^M$. An analogue of the Koebe function for this class is the
Pick function
$$P_M(z)=MK^{-1}(K(z)/M)=z+\sum_{n=1}^{\infty}p_n(M)z^n.$$ The
sharp estimate $|a_2|\leq 2(1-1/M)=p_2(M)$ in the class $S^M$ is
rather trivial and has been obtained by G.~Pick [12] in 1917. The
next coefficient $a_3$ was estimated independently by
A.~C.~Schaeffer and D.~C.~Spencer [14] in 1945 and O.~Tammi [16]
in 1953. The Pick function does not give the maximum to $|a_3|$
and the estimate is much more difficult. M.~Schiffer and O.~Tammi
[15] in 1965 found that $|a_4|\leq p_4(M)$ for any $f\in S^M$ with
$M>300$. This result was repeated by O.~Tammi [18, page 210] in a
weaker form ($M>700$) and there it was conjectured that this
constant could be decreased until 11. The case of function with
real coefficients is simpler: the Pick function gives the maximum
to $|a_4|$ for $M\geq 11$ and this constant is sharp (see [17],
[19, p.163]). By our suggested method we will show that the Pick
function locally maximizes $|a_4|$ on $S^M$ if $M>
M_0=22.9569\dots$ and does not for $1< M<M_0$. This disproves
Tammi's conjecture.

\head{\S 2. Preliminary statements}\endhead

The parametric representation of univalent functions is based on
the L\"owner differential equation and goes back to the famous
L\"owner's paper [11] where the author using an idea of semigroups
of conformal maps derived the equation $$
\frac{dw}{dt}=-w\frac{e^{iu}+w}{e^{iu}-w},\;\;\;w|_{t=0}=z,\;\;\;t\geq0,
$$ where the control function $u=u(t)$ is piecewise continuous in
$t\geq0$. One finds the foundations of the parametric method,
e.g., in [1,8,13]. It is convenient to make the change of
variables $t\to 1-e^{-t}$ and rewrite L\"owner's equation using
the preceding notation for the independent variable as follows $$
\frac{dw}{dt}=\frac{-w}{(1-t)}\frac{e^{iu}+w}{e^{iu}-w},\;\;\;w|_{t=0}=z,
\;\;\;0\leq t\leq1.\tag1 $$ A.~C.~Schaeffer and D.~C.~Spencer [14]
were the first who used L\"owner's equation for the class $S^M$
and proved that the integrals $$
w=w(z,t)=(1-t)(z+a_2(t)z^2+\dots)\tag2 $$ of the equation (1)
represent a dense subclass of functions $f\in S^M$ by $$
f(z)=Mw(z,1-1/M).\tag3 $$  Representation (3) is valid for all
$M\geq1$ including $M=\infty$ if the product in (3) is regarded as
the limit as $M\to\infty$. From now on, we will use the notation
(2).

We remark that the case $u=\pi$ in (1) corresponds to the Koebe
function in the class $S$ or to the Pick function in the class
$S^M$ by (3). Besides, the dense subclass of $S^M$ represented by
(3) contains all functions that give the boundary points of the
coefficient region $$ V_4^M=\{(a_2,a_3,\R a_4):f\in
S^M\},\;\;\;1\leq M\leq\infty. $$

For given real numbers $\mu$ and $\nu$, we will consider the
linear functional $$ L(\mu,\nu;f)=a_2+\mu a_3+\nu a_4 $$ in $S^M$.
The Koebe function $K(z)$ maximizes $\R L(0,0;f)$ in $S$, and
similarly, the Pick function $P_M(z)$ maximizes $\R L(0,0;f)$ in
$S^M$. We will describe the set of  $\mu$ and $\nu$ for which the
local maximum of $\R L(\mu,\nu;f)$ in $S^M$ is attained by the
Pick function and  will apply this result to the extremal problems
for the class $S$ or $S^M$ stated in the introduction.

We write $a_k(t)\equiv x_{2k-3}(t)+ix_{2k-2}(t)$, $k=2,3,4$.
Substituting (2) into (1) we obtain the following differential
equations $$\align
 \dot x_1(t) &=-2\cos u,\;\;\;x_1(0)=0, \\
  \dot x_2(t) &=2\sin u,\;\;\;x_2(0)=0, \\
   \dot x_3(t) &=-4(x_1\cos
u+x_2\sin u)+2(t-1)\cos2u,\;\;\;x_3(0)=0,\tag 4\\
 \dot x_4(t) &=4(x_1\sin
u-x_2\cos u)-2(t-1)\sin2u,\;\;\;x_4(0)=0, \\
 \dot x_5(t) &=-2((2x_3+x_1^2-x_2^2)\cos u+2(x_4+x_1x_2)\sin u)\\
 &+ 6(t-1)
(x_1\cos2u+x_2\sin2u)-2(t-1)^2\cos3u,\;\;\;  x_5(0)=0. \endalign$$
The extremal problem $$ \R L(\mu,\nu;f)\to\max $$ in the class
$S^M$ is equivalent to the extremal problem $$ x_1(1-1/M)+\mu
x_3(1-1/M)+\nu x_5(1-1/M)\to\max\tag5 $$ for solutions to the
system (4). The parametric representation (4) for the coefficients
generated by  L\"owner's equation allows us to apply the classical
variational methods [3] or  Pontryagin's maximum principle [10].
We introduce the Hamiltonian function in order to formulate the
necessary extremum conditions for the problem (5) $$\align
H(t,x,\Psi,u)&=-2\cos u\Psi_1+2\sin u\Psi_2\\ &-(4(x_1\cos
u+x_2\sin u)-2(t-1)\cos2u)\Psi_3\\ &+(4(x_1\sin u-x_2\cos
u)-2(t-1)\sin2u)\Psi_4\tag 6\\ &-(2((2x_3+x_1^2-x_2^2)\cos
u+2(x_4+x_1x_2)\sin u)\\
&-6(t-1)(x_1\cos2u+x_2\sin2u)+2(t-1)^2\cos3u)\Psi_5,
\endalign
$$
where $x=(x_1,\dots,x_5)^T$ satisfies (4) and
$\Psi=(\Psi_1,\dots,\Psi_5)^T$ satisfies the conjugate system
$$\aligned
\dot\Psi_1&=4\cos u\Psi_3-4\sin u\Psi_4+(4x_1\cos u+4x_2\sin u-
6(t-1)\cos2u)\Psi_5,\\
\dot\Psi_2&=4\sin u\Psi_3+4\cos u\Psi_4-(4x_2\cos u-4x_1\sin u+
6(t-1)\sin2u)\Psi_5,\\
\dot\Psi_3&=4\cos u\Psi_5,\\
\dot\Psi_4&=4\sin u\Psi_5,\\
\dot\Psi_5&=0,
\endaligned\tag7
$$
and the transversality conditions
$$\aligned
\Psi_1(1-1/M)&=1,\;\;\Psi_3(1-1/M)=\mu,\;\;\Psi_5(1-1/M)=\nu,\\
\Psi_2(1-1/M)&=\Psi_4(1-1/M)=0.\endaligned\tag8
$$

The optimal control function $u^*$ corresponding to the extremal
function $f^*$ in (5) satisfies Pontryagin's maximum principle
$$
\max_uH(t,x^*,\Psi^*,u)=H(t,x^*,\Psi^*,u^*),\;\;\;0\leq t\leq1-1/M,\tag9
$$
where $(x^*,\Psi^*)$ is the solution to (4) and (7) with $u=u^*$ in their
right-hand sides. Hence, $u^*$ is a root of the equation
$$
H_u(t,x,\Psi,u)=0\tag10
$$
for $x=x^*$ and $\Psi=\Psi^*$.

\proclaim{Lemma 1} Let us suppose that a control function $u$ in
(4) and (7-8) generates the solutions $x(t)$ and $\Psi(t)$ for
which $u$ satisfies (9), is unique up to the $2\pi$-translation,
and $$ H_{uu}(t,x,\Psi,\pi)\neq0,\;\;\;0\leq t\leq1-1/M.\tag11 $$
Let us denote by $(x(t,\xi),\Psi(t,\xi))$  solutions to (4) and
(7) with the initial conditions $\Psi(0,\xi)=\Psi(0)+\xi$ and
$u=u(t,\xi)$ in their right-hand sides satisfying the maximum
principle (9). Then for $\xi\to0$, we have the following
asymptotic behaviour $$
\|(x(1-1/M,\xi),\Psi(1-1/M,\xi))-(x(1-1/M),\Psi(1-1/M))\|=o(1), $$
where $\|\cdot\|$ is the Euclidean vector norm.
\endproclaim
\demo{Proof}
 Since there exists a unique solution $u$ satisfying (9) and
(11), the same is true for a slightly changed parameters of the
function $H$. Therefore, equations (9) and (10) locally determine
a unique continuous implicit function $u=u(t,x,\Psi)$ satisfying
the maximum principle. Writing
$u(t,\xi)=u(t,x(t,\xi),\Psi(t,\xi))$ we substitute it into (4) and
(7). Now we apply the theorem on the continuous dependence of
solutions of differential equations on the initial conditions and
complete the proof of Lemma 1. \quad\qed\enddemo

If the Pick function $P_M$ is extremal for (5), then $u=\pi$ is the optimal
control function, (4) and (7) give
$(x(t),\Psi(t))=(x^0(t),\Psi^0(t))$, where
$$
x_1^0(t)=2t,\;\;\;x_3^0(t)=5t^2-2t,\;\;\;x_2^0(t)=x_4^0(t)=0,\tag12
$$
and
$$\aligned
\Psi_1^0(t)&=\nu\left(t-1+\frac{1}{M}\right)^2+\left(\frac{14\nu}{M}-8\nu-
4\mu\right)\left(t-1+\frac{1}{M}\right)+1,\\
\Psi_3^0(t)&=-4\nu\left(t-1+\frac{1}{M}\right)+\mu,\;\;\;\Psi_5^0(t)=\nu,
\;\;\;\Psi_2^0(t)=\Psi_4^0(t)=0.
\endaligned\tag13$$
The conditions of Lemma 1 play the key role as a necessary local
extremum condition for the Pick function $P_M$. To verify these we
substitute (12) and (13) into the Hamiltonian function
$H(t,x,\Psi,u)$ given by (6) and study the extremum properties of
$H(t,u)=H(t,x^0(t),\Psi^0(t),u)$ which is just a cubic polynomial
of $\cos u$. Let us describe a set of suitable real parameters
$(\mu,\nu)$ satisfying Lemma 1. Let $D(M)$ denote the maximal
domain in the $(\mu,\nu)$-plane which is starlike with the respect
to the origin and satisfies the following conditions: \roster
\item"[i]" $H(t,u)$ as a function of $y=\cos u$ attains its
maximum on $[-1,1]$ only at $y=-1$ for all $t\in[0,1-1/M]$;
\item"[ii]" $H_{uu}(t,\pi)\neq0$, $0\leq t\leq1-1/M$.
\endroster
We will consider $(\mu,\nu)\in D(M)$. The point $x^0(1-1/M)$
belongs to the boundary $\partial V_4^M$ of $V_4^M$ and is given
by the Pick function. Each $x\in\partial V_4^M$ can be obtained as
a solution to the system (4) with a certain optimal control
function $u$. Lemma 1 admits a reverse formulation.

\proclaim{Lemma 2} Let $(\mu,\nu)\in D(M)$ and denote by
$(x^-(t,\xi),\Psi^-(t,\xi))$  solutions to (4) and (7) with the boundary
conditions $\Psi^-(1-1/M,\xi)=\Psi^0(1-1/M)+\xi$. Let $u=u^-(t,\xi)$ in
their right-hand sides satisfy the maximum principle (9). If
$$
\|x^-(1-1/M,\xi)-x^0(1-1/M)\|=o(1),\;\;\;\text{as $\xi\to0$},
$$
then
$$
\|\Psi^-(0,\xi)-\Psi^0(0)\|=o(1).
$$
\endproclaim
The {\it proof} of Lemma 2 is similar to that of Lemma 1 reversing the direction of
variation of $t$ from $1-1/M$ to 0 and noting that $x^-(0)=x^0(0)=0$.

Lemmas 1 and 2 imply that if $(x(t),\Psi(t))$ is given by (4) and (7),
$(x(1-1/M),\Psi(1-1/M))$ is close to $(x^0(1-1/M),\Psi^0(1-1/M))$, and
$x(1-1/M)\in\partial V_4^M$, then $(x(t),\Psi(t))$ is equal to
$(x(t,\xi),\Psi(t,\xi))$ for a certain $\xi$ close to 0.

The principles of calculus of variations interpret geometrically the
transversality conditions as an orthogonality property of $\Psi(1-1/M)$
to all possible variations of $x(1-1/M)$ in $\partial V_4^M$.
Nevertheless, we rigorously prove this fact for completeness.

\proclaim{Lemma 3} Let us suppose that $(\mu,\nu)\in D(M)$ and the
initial conditions in (7) are
$$
\Psi(0,\xi)=\Psi^0(0)+\xi,\;\;\;\xi=\varepsilon e,\;\;\;\varepsilon>0,\;\;\;
e=(e_1,\dots,e_5)^T,\;\;\;\|e\|=1,\tag14
$$
and
$$
x(t,\xi)=x^0(t)+\varepsilon\delta x(t)+o(\varepsilon),\;\;\;\varepsilon\to0.\tag15
$$
Then $\Psi^0(1-1/M)$ is orthogonal to $\delta x(1-1/M)$, if
$\delta x(1-1/M)\neq0$.
\endproclaim
\demo{Proof}
First we note that the conditions of Lemma 1 guarantee
the differentiability of $x(t,\varepsilon e)$ with respect to $\varepsilon$ at
$\varepsilon=0$. Therefore, the representation (15) is valid and according to
Lemmas 1 and 2 the expansion (15) produces all possible variations
$\delta x(1-1/M)$ at $x^0(1-1/M)$ associated with $\Psi^0(1-1/M)$.

We denote the column of the functions in the right-hand side of
(4) by $g(t,x,u)$ and rewrite (4) in the vector form as $$ \dot
x=g(t,x,u),\;\;\;x(0)=0.\tag16 $$ The system (7) is equivalent to
$$ \dot\Psi=-\frac{\partial H}{\partial
x}(t,x,\Psi,u),\;\;\;\Psi(0)= \Psi^0(0)+\varepsilon e^T.\tag17 $$
Substituting (15) into (16) we obtain $$ \frac{d{\delta
x}}{dt}=g_x\delta x+g_uu_{\varepsilon}, $$ by differentiating with
respect to $\varepsilon$ at $\varepsilon=0$. This, together with
(4), (6), (10), and (17), imply that $$ \frac{d}{dt}((\delta
x)^T\Psi)=(\delta x)^Tg_x^T\Psi+ u_{\varepsilon}g_u^T\Psi+(\delta
x)^T\dot\Psi=0, $$ because the second term is equal to
$u_{\varepsilon}H_u(t,x,\Psi,u)=0$, and the remaining terms give
the zero sum since $\dot\Psi=-g_x^T\Psi$.

Thus, we see that $(\delta x)^T\Psi$ does not depend on $t$ and
vanishes at $t=0$ because $\delta x(0)=0$. This completes the proof of
Lemma 3.
\quad\qed\enddemo

\head \S 3. Local extremum conditions\endhead

\proclaim{Lemma 4} Under the conditions of Lemma 3 we suppose that
in the vector $e=(e_1,\dots,e_5)^T$, which corresponds to the
variation of $\Psi^0(0)$ in (14), the coordinates $e_2$ and $e_4$
vanish. Then $\delta x(1-1/M)=0$.\endproclaim \demo{Proof} The
condition $e_2=e_4=0$ implies that the systems (4) and (7) have
vanishing coordinate solutions $x_2(t)=x_4(t)=0$ and
$\Psi_2(t)=\Psi_4(t)=0$. In this case the Hamilton function $H$ is
a polynomial of $y=\cos u$, which has a unique maximum on $[-1,1]$
at $y=-1$. Its derivative with respect to $y$ does not vanish  at
$y=-1$ for $\varepsilon>0$ sufficiently small. Hence, $u\equiv\pi$
is a unique optimal control function for such $\varepsilon$ and
$x(1-1/M,\xi)=x^0(1-1/M)$, that ends the proof of Lemma 4.
\quad\qed\enddemo

By  analogy with the expansion (15) in Lemma 3
we have the expansion
$$
\Psi(t,\xi)=\Psi^0(t)+\varepsilon\delta\Psi(t)+o(\varepsilon),\;\;\;\varepsilon\to0.
$$
 Lemma 4 shows that the condition $e_2=e_4=0$ implies $\delta x=0$
and $\delta\Psi_2=\delta\Psi_4=0$. Only $\Psi_1$, $\Psi_3$,
and $\Psi_5$ can vary in this case . It follows  from Lemmas 1-4 that we should consider
variations $\Psi(0,\xi)$ by (14) with $e_1=e_3=e_5=0$ in order to study the
character of the point $x^0(1-1/M)\in\partial V_4^M$.

Let us set $\xi=(0,p,0,q,0)^T$ with arbitrary real $p$ and $q$,
and study $x(1-1/M,\xi)$ in a neighborhood of $x^0(1-1/M)$. In
other words, we will solve the systems (4) and (7) with the
initial conditions $x(0)=0$ and (14), which we rewrite with
coordinates
$$\aligned
\Psi_1(0)&=3\nu\left(1-\frac{1}{M}\right)\left(3-\frac{5}{M}\right)+
4\mu\left(1-\frac{1}{M}\right)+1,\\
\Psi_2(0)&=p,\\
\Psi_3(0)&=4\nu\left(1-\frac{1}{M}\right)+\mu,\\
\Psi_4(0)&=q,\\
\Psi_5(0)&=\nu.
\endaligned\tag 18$$
Let
$$
F:(p,q)\to x_1(1-1/M)+\mu x_3(1-1/M)+\nu x_5(1-1/M)
$$
be a real valued mapping from the $(p,q)$-plane onto the linear
combination of the components of the solution to the Cauchy
problem for the systems (4) and (7) with the initial conditions
(18). The control function $u$ in the right-hand side of (4) and
(7) satisfies the maximum principle. The mapping $F$ is well
defined in a neighborhood of $(0,0)$ if $(\mu,\nu)\in D(M)$. In
this case $u=u(t,x,\Psi)$ is an implicit function defined by (10).
Since $(x,\Psi)$ depends only on $(p,q)$, we denote by
$u(t,p,q)=u(t,x(p,q),\Psi(p,q))$. We note that $F(0,0)=\R
L(\mu,\nu;P_M)$ and the values of $F(p,q)$ correspond to those of
$\R L(\mu,\nu;f)$ with respect to the variations of $P_M$
generated by the initial conditions $\Psi_2(0)=p$ and
$\Psi_4(0)=q$.

\proclaim{Theorem 5} Let us suppose that $(\mu,\nu)\in D(M)$. If $P_M$ locally
maximizes $\R L(\mu,\nu;f)$ in $S(M)$, then
$$
F_{pp}(0,0)\leq0,\;\;\;F_{pp}(0,0)F_{qq}(0,0)-F_{pq}^2(0,0)\geq0.
$$
Conversely, if
$$
F_{pp}(0,0)<0,\;\;\;F_{pp}(0,0)F_{qq}(0,0)-F_{pq}^2(0,0)>0,
$$
then $P_M$ locally maximizes $\R L(\mu,\nu;f)$ in $S(M)$.
\endproclaim
\demo{Proof} We first claim that $(0,0)$ is a critical point of
$F(p,q)$. Indeed, substituting $u=u(t,p,q)$ in the three equations
in (4) (for $x_1,x_3,x_5$) and differentiating them with respect
to $p$ and $q$, we obtain differential equations for $(x_k)_p$ and
$(x_k)_q$, $k=1,3,5$, with vanishing initial conditions.
Substituting there $p=q=0$ and $u=\pi$, $x_2(t)=x_4(t)=0$ we find
that all derivatives $(\dot x_k)_p$ and $(\dot x_k)_q$, $k=1,3,5$
are identically zeros and, hence, $F_p(0,0)=F_q(0,0)=0$.

The first statement of Theorem 5 means that the quadratic form of
the second differential of $F$ at $(0,0)$ is negatively
semi-definite which is the necessary  condition of local extremum.
Similarly, the second statement signifies that the above quadratic
form is negative definite which is the sufficient local extremum
condition. This completes the proof of Theorem 5.
\quad\qed\enddemo

The same reasoning can be made for the class $S_R^M$ of functions
$f\in S^M$ with real Taylor coefficients $a_n$, $n\geq2$. The
coefficients of an arbitrary boundary function $f_R\in S_R^M$ for
the set $V_{R4}^M=\{(a_2,a_3,a_4):f\in S_R^M\}$ can be obtained by
integrating the systems (4) and (7) with the control function $u$
satisfying the maximum principle (9) and with vanishing initial
values of $\Psi_2(0)$ and $\Psi_4(0)$. Therefore, variations of
$\Psi_2(0)$ and $\Psi_4(0)$ are forbidden and the Pick function
$P_M\in S_R^M$ locally maximizes $L(\mu,\nu;f)$ in $S_R^M$ if
$(\mu,\nu)\in D(M)$.

Now we will apply Theorem 5 to construct an analytic and numerical solution
process.
We need to calculate the
partial derivatives $u_p$ and $u_q$ at $(0,0)$ to evaluate the
 partial derivatives of $F$ at $(0,0)$. Differentiating (10) with
respect to $p$ and $q$ we obtain $$\align
H_{ux}x_p+H_{u\Psi}\Psi_p+H_{uu}u_p&=0,\\
H_{ux}x_q+H_{u\Psi}\Psi_q+H_{uu}u_q&=0,\endalign $$ which leads us
to the formulae $$
u_p=-\frac{H_{ux}x_p+H_{u\Psi}\Psi_p}{H_{uu}}\tag19 $$ and $$
u_q=-\frac{H_{ux}x_q+H_{u\Psi}\Psi_q}{H_{uu}}.\tag20 $$ Direct
calculation gives $$ H_{uu}(t,x^0,\Psi^0,\pi)=-2\left[16\nu
t^2-4\left(2\nu+\frac{4\nu}{M}-
\mu\right)t+2\nu+1-\frac{4(2\nu+\mu)}{M}+\frac{15\nu}{M^2}\right].\tag21
$$ Differentiating (6) with respect to corresponding variables at
$u=\pi$ we find $$\align
H_{ux_2}(t,x^0,\Psi^0,\pi)&=4\left(\nu\left(t+1-\frac{4}{M}\right)+
\mu\right),\tag22\\ H_{ux_4}(t,x^0,\Psi^0,\pi)&=4\nu,\tag23\\
H_{u\Psi_2}(t,x^0,\Psi^0,\pi)&=-2,\tag24\\
H_{u\Psi_4}(t,x^0,\Psi^0,\pi)&=4(1-3t).\tag25
\endalign$$
Differentiating (7) with respect to $p$ and $q$ at $(0,0)$ with
$u=\pi$, $(x_1)_p(t)=(x_1)_q(t)=x_2(t)=0$ we see that
$$
(\Psi_1)_p(t)=(\Psi_1)_q(t)=(\Psi_3)_p(t)=(\Psi_3)_q(t)=(\Psi_5)_p(t)=
(\Psi_5)_q(t)=0.\tag26
$$
The formulae (22--26) allow us to calculate the numerators in (19) and (20) as
$$
4\left(\nu\left(t+1-\frac{4}{M}\right)+\mu\right)(x_2)_p+4\nu(x_4)_p-
2(\Psi_2)_p+4(1-3t)(\Psi_4)_p
$$
and
$$
4\left(\nu\left(t+1-\frac{4}{M}\right)+\mu\right)(x_2)_q+4\nu(x_4)_q-
2(\Psi_2)_q+4(1-3t)(\Psi_4)_q
$$
respectively.

From (4) and (7) we conclude that $\dot\Psi_4=2\nu\dot x_2$ yields the
equalities ${(\dot\Psi_4)}_p=2\nu{(\dot x_2)}_p$ and
${(\dot\Psi_4)}_q=2\nu{(\dot x_2)}_q$. The initial conditions $(\Psi_4)_p(0)=0$ and
$(\Psi_4)_q(0)=1$ imply that $(\Psi_4)_p=2\nu(x_2)_p$
and $(\Psi_4)_q=2\nu(x_2)_q+1$.

We substitute the last relations and (21) into (19--20), and finally, get
$$
u_p=\frac{((3-5t-4/M)\nu+\mu)2y_4+2\nu y_5-y_6}{16\nu t^2-4(2\nu+
4\nu/M-\mu)t+2\nu+1-4(2\nu+\mu)/M+15\nu/M^2},\tag27
$$
where $y_4:=(x_2)_p$, $y_5:=(x_4)_p$, $y_6:=(\Psi_2)_p$, and
$$
u_q=\frac{((3-5t-4/M)\nu+\mu)2y_{10}+2\nu y_{11}-y_{12}+2(1-3t)}
{16\nu t^2-4(2\nu+4\nu/M-\mu)t+2\nu+1-4(2\nu+\mu)/M+15\nu/M^2},\tag28
$$
where $y_{10}:=(x_2)_q$, $y_{11}:=(x_4)_q$, $y_{12}:=(\Psi_2)_q$.

Set $y_1:=(x_1)_{pp}$, $y_2:=(x_3)_{pp}$, and $y_3:=(x_5)_{pp}$.
Differentiating (4) twice with respect to $p$ at $u=\pi$,
$x_2=x_4=\Psi_2=\Psi_4=(x_1)_p=(x_3)_p=0$ we obtain
$$\align
\dot y_1&=-2u_p^2,\;\;\;y_1(0)=0,\tag29\\
\dot y_2&=4(y_1+2y_4u_p-2(2t-1)u_p^2),\;\;\;y_2(0)=0,\tag30\\
\dot y_3&=2(7t-3)y_1+4y_2-4y_4^2+8(5t-3)y_4u_p+8y_5u_p\\&-
2(47t^2-46t+9)u_p^2,\;\;\;y_3(0)=0.\tag31
\endalign
$$
Similarly, we differentiate the remaining equations in (4) with respect to $p$
and obtain
$$\align
\dot y_4&=-2u_p,\;\;\;y_4(0)=0,\tag32\\
\dot y_5&=4(y_4+(1-3t)u_p),\;\;\;y_5(0)=0.\tag33
\endalign
$$
Finally, we differentiate the second equation in (7) with respect to $p$ and
get
$$
\dot y_6=-4\nu\left(\left(t+1-\frac{4}{M}+\mu/\nu\right)u_p+y_4\right),
\;\;\;y_6(0)=1.\tag34
$$
Summarizing, we have deduced an evaluation algorithm for $F_{pp}(0,0)$ expressed
by the following theorem.

\proclaim{Theorem 6} Suppose $(\mu,\nu)\in D(M)$. Let
$y_1(t),\dots,y_6(t)$, $0\leq t\leq1-1/M$, be solutions to the
Cauchy problem  for the differential equations (29--34). Then the
relation $$ F_{pp}(0,0)=y_1(1-1/M)+\mu y_2(1-1/M)+\nu y_3(1-1/M)
$$ is valid.
\endproclaim

\remark{Remark} The subsystem (32--34) can be solved independently because these
equations do not contain $y_1$, $y_2$, and $y_3$.
\endremark

Calculation of $F_{qq}(0,0)$ and $F_{pq}(0,0)$ may be handled in
much the same way. Let $y_7:=(x_1)_{qq}$, $y_8:=(x_3)_{qq}$, and
$y_9:=(x_5)_{qq}$. From (4) we have $$\align \dot
y_7&=-2u_q^2,\;\;\;y_7(0)=0,\tag35\\ \dot
y_8&=4(y_7+2y_{10}u_q-2(2t-1)u_q^2),\;\;\;y_8(0)=0,\tag36\\ \dot
y_9&=2(7t-3)y_7+4y_8-4y_{10}^2\\ &+8(5t-3)y_{10}u_q+8y_{11}u_q-
2(47t^2-46t+9)u_q^2,\;\;\;y_9(0)=0.\tag37\endalign $$
Differentiating the two even equations in (4) and the second
equation in (7) with respect to $q$ we get $$\align \dot
y_{10}&=-2u_q,\;\;\;y_{10}(0)=0,\tag38\\ \dot
y_{11}&=4(y_{10}+(1-3t)u_q),\;\;\;y_{11}(0)=0,\tag39\\ \dot
y_{12}&=-4\nu\left(\left(t+1-\frac{4}{M}+\mu\right)u_q+y_{10}\right)-4,
\;\;\;y_{12}(0)=0.\tag40\endalign $$ Let $y_{13}:=(x_1)_{pq}$,
$y_{14}:=(x_3)_{pq}$, and $y_{15}:=(x_5)_{pq}$. We continue in
this fashion differentiating (4) with respect to $p$, and
subsequently, with respect to $q$ and obtain $$\align \dot
y_{13}&=-2u_pu_q,\;\;\;y_{13}(0)=0,\tag41\\ \dot
y_{14}&=4(y_{13}+y_4u_q+y_{10}u_p)-8(2t-1)u_pu_q,\;\;\;y_{14}(0)=0,
\tag42\\ \dot
y_{15}&=2(7t-3)y_{13}+4y_{14}-4y_4y_{10}+4(5t-3)(y_4u_q+y_{10}u_p)\\
&+4y_5u_q+4y_{11}u_p-2(47t^2-46t+9)u_pu_q,\;\;\;y_{15}(0)=0.\tag43\endalign
$$

Summing up the calculation process we derive the following theorem.

\proclaim{Theorem 7} Suppose $(\mu,\nu)\in D(M)$. Let $
y_7(t),\dots,y_{12}(t)$, $0\leq t\leq1-1/M$, be solutions to the
Cauchy problem for the differential equations (35--40). Then, the
relation $$ F_{qq}(0,0)=y_7(1-1/M)+\mu y_8(1-1/M)+\nu y_9(1-1/M)
$$ holds. Let $y_4(t), y_5(t), y_6(t)$, and $y_{13}(t), y_{14}(t),
y_{15}(t)$, $0\leq t\leq1-1/M$, be solutions to the Cauchy problem
for the differential equations (32--34) and (41--43) respectively.
Then, the relation $$ F_{pq}(0,0)=y_{13}(1-1/M)+\mu
y_{14}(1-1/M)+\nu y_{15}(1-1/M) $$ holds.
\endproclaim

\remark{Remark} As in the remark after Theorem 6 we note that the subsystem (38--40) can be solved independently because these
equations do not contain $y_7$, $y_8$, and $y_9$.
\endremark

\head \S 4. Explicit integration\endhead

An explicit integration of the systems in Theorems 6 and 7 is
possible only in the case $\nu=0$, i.e., when the linear
functional $L$ does not depend on $a_4$. In this case the two last
equations in (4) and (7) disappear and the mapping $F$ becomes a
function of a variable $p$ as
$$
F:p\to x_1(1-1/M)+\mu x_3(1-1/M).
$$
The criterion in Theorem 5 is reduced to the
inequality $F''(0)<0$, where $F''(0)=y_1(1-1/M)+\mu y_2(1-1/M)$. The systems
(29--34), (35--40), and (41--43) are reduced to the four equations (29),
(30), (32), and (34). We substitute $\nu=0$ into (34, 27) and obtain
that
$$
\dot y_6=-4\mu u_p,\;\;\;y_6(0)=1,\tag44
$$
and
$$
u_p=\frac{2\mu y_4-y_6}{4\mu t+1-4\mu/M}.\tag45
$$
Compare (32, 44) and observe that $\dot y_6=2\mu\dot y_4$, which
implies
$$
y_6(t)=2\mu y_4(t)+1.
$$
This equation allows us to exclude $y_6$ from (45) and integrate the
system of equations (29), (30), and (32) with
$$
u_p=\frac{-1}{4\mu t+1-4\mu/M}.\tag46
$$
The equations (29) and (32) give
$$
y_1(t)=\frac{1}{2\mu}\left(\frac{1}{1-4\mu/M+4\mu t}-\frac{1}
{1-4\mu/M}\right)\tag47
$$
and
$$
y_4(t)=\frac{1}{2\mu}\log\frac{1-4\mu/M+4\mu t}{1-4\mu/M}.\tag48
$$
Substituting (46) and (48) into (30) we integrate it, then, taking into account
(47), we finally obtain that
$$
y_1(t)+\mu y_2(t)=\frac{-1}{2\mu}\left[\log^2((1+4\mu-8\mu/M)(1-4\mu/M))+
\log\frac{1+4\mu-8\mu/M}{1-4\mu/M}\right].\tag49
$$
Making use of (49) we formulate the following theorem as a corollary of Theorem 6.

\proclaim{Theorem 8} Suppose $(\mu,0)\in D(M)$ and $\mu$ satisfies the
inequality
$$
\frac{1}{2\mu}\left[\log^2((1+4\mu-8\mu/M)(1-4\mu/M))+
\log\frac{1+4\mu-8\mu/M}{1-4\mu/M}\right]>0.\tag50
$$
Then, the Pick function $P_M$ locally maximizes
$\R L(\mu,0;f)=\R (a_2+\mu a_3)$ in $S^M$.
If the left-hand side of (50) is negative, then $P_M$ does not give a
local maximum to $\R L(\mu,0;f)$ in $S^M$.
\endproclaim

We note that $\{\mu:(\mu,0)\in D=D(\infty)\}=\{\mu:\mu>-0.25\}$. If
$M=\infty$,
then the inequality (50) is of the form
$$
\frac{-1}{2\mu}(1+\log(1+4\mu))\log(1+4\mu)<0,
$$
that is equivalent to the best possible inequality
$\mu>-\lambda_0:=-(e-1)/{4e}$. This means that the Bombieri number
$\sigma_{32}$ is equal to $\lambda_0$. This result has been
obtained recently by R.~Greiner and O.~Roth [9].

\head \S 5. Bombieri's number $\sigma_{42}$\endhead

Now we apply Theorems 5-7 to evaluate Bombieri's number
$\sigma_{42}$.

\proclaim{Proposition} Let $m,n\geq 2$ be fixed integers. Then
$$
\sigma_{mn}=\sup\{\lambda\in \Bbb R:\,\,\R(a_n-\lambda a_m),
\text{is locally maximized on $S$ by K(z)}\}.
$$
\endproclaim

The proof of this statement is quite obvious and can be found,
e.g., in [9].

 According to Theorem 5, Bombieri's number
$\sigma_{42}$ is calculated as
$$
-\inf\{\nu':F_{pp}(0,0)<0\;\text{and}\;F_{pp}(0,0)F_{qq}(0,0)-
F_{pq}^2(0,0)>0\;\text{for}\;\mu=0,\;M=\infty,\;\nu\in[\nu',0]\}.
$$
Theorems 6 and 7 reduce the problem to the solution of the equations
$$
y_1(1)+\nu y_3(1)=0,\tag51
$$
or
$$
(y_1(1)+\nu y_3(1))(y_7(1)+\nu y_9(1))-(y_{13}(1)+\nu y_{15}(1))^2=0,\tag52
$$
where $y_1(t),\dots,y_{15}(t)$ are solutions to the Cauchy problem
for the differential equations (29--43) with $\mu=0$ and
$M=\infty$. Thus we are able to formulate the following result.

\proclaim{Theorem 9} Bombieri's number $\sigma_{42}$ is equal to
 the maximum of the negative roots of the equations (51) and
(52) multiplied by (-1), where $y_1(t),\dots,y_{15}(t)$ are
solutions to (29--43) for $\mu=0$ and $M=\infty$.
\endproclaim

To illustrate this result we give the numerical approximation
$\sigma_{42}\approx 0.050057...$

\remark{Remark} We  used Wolfram's {\it Mathematica} to evaluate
numerically $\sigma_{42}$ as well as other  Bombieri's numbers in
the next sections.  More precisely, the combination of 4-th order
Runge-Kutta and Adams methods is used. It allows to reach a higher
precision at short computational time. Actually, the level of
precision imply the number of iterations and can be prescribed as
the machine precision.  For solving the above systems the level
$10^{-12}$ is reached in time of order of few minutes on a machine
with a Pentium 4 processor (1800~Mhz) and having 512~Mb of RAM.
\endremark

\remark{Remark} Note that $\{\nu:(0,\nu)\in
D=D(\infty)\}=\{\nu:\nu>-0.1\}$. O.~Roth showed (private
communication) that $\sigma_{42}\leq0.050284...$, the latter
number is equal to $(2-\R a_2)/(4-\R a_4)$ for a function sequence
which is critical in the problem of finding $\sigma_{32}$.
\endremark

\head \S 6. Bombieri's numbers $\sigma_{24}$ and $\sigma_{34}$
\endhead

The problem of finding Bombieri's number $\sigma_{24}$ is reduced to
description of $\nu$ for which the Koebe function $K$ minimizes
$\R L(0,\nu;f)$ in $S$. To follow the preceding scheme we will write $L^1(\nu;f)=a_4+\nu a_2$
in place of $L$. Now $\sigma_{24}$ is equal to
the supremum over all real values $\lambda_{24}$, such that
$\R L^1(-\lambda_{24};f)$ is locally maximized by the Koebe function $K$
in $S$. Again we consider the system of equations (4), the Hamilton function
$H(t,x,\Psi,u)$ given by (6), and the system of equations (7), where the
transversality conditions (8) are replaced by
$$
\Psi_1(1)=\nu,\;\;\;\Psi_5(1)=1,\;\;\;\Psi_2(1)=\Psi_3(1)=\Psi_4(1)=0.\tag53
$$
For $u=\pi$ let us denote the integrals for (4) and (7) by $(x^0(t),\Psi^1(t))$, where
$x_1^0(t),\dots,x_4^0(t)$ are given by (12) and
$$
\Psi_1^1(t)=t^2-10t+9+\nu,\;\;\;\Psi_3^1(t)=-4(t-1),\;\;\;\Psi_5^1(t)=1,
\;\;\;\Psi_2^1(t)=\Psi_4^1(t)=0.\tag54
$$
As before, the control function $u=u(t)$ satisfies the maximum principle
(9) and, hence, the equation (10). We consider only $\nu\in D^1$, where $D^1$
denotes a maximal interval on the $\nu$-axis that satisfies the following
conditions:
\roster
\item"[i]" $H^1(t,u)=H(t,x^0(t),\Psi^1(t),u)$ as a function of $y=\cos u$ attains
its maximum on $[-1,1]$ only at $y=-1$ for all $t\in[0,1]$;
\item"[ii]"  $H_{uu}^1(t,\pi)\neq0$, $0\leq t\leq1$.
\endroster
Note that $D^1=\{\nu:\nu>-1\}$.

Similarly to Section 3 we solve the systems (4) and (7) with the initial
conditions $x(0)=0$ and
$$
\Psi_1(0)=9+\nu,\;\;\;\Psi_2(0)=p,\;\;\;\Psi_3(0)=4,\;\;\;\Psi_4(0)=q,
\;\;\;\Psi_5(0)=1,\tag55
$$
where $p$ and $q$ are arbitrary real numbers. The initial conditions (55) are thought of as
variations of (54) at $t=0$.

Let
$$
F^1:(p,q)\to\nu x_1(1)+x_3(1)
$$
be a real valued mapping from the $(p,q)$-plane to the linear combination of
the components of a solution to the Cauchy problem for (4) and (7) with
the initial conditions (55) and the control $u$ satisfying the maximum principle.
The mapping $F^1$ is well defined in a neighborhood of $(0,0)$ if
$\nu\in D^1$, and $u=u(t,x,\Psi)$ is an implicit function defined by
(10). We denote by $u(t,p,q)=u(t,x(p,q),\Psi(p,q))$.
Lemmas 1-4 and Theorem 5 admit an analogous formulation
for the functional $L^1(\nu;f)$. So,
if $K$ locally maximizes $\R L^1(\nu;f)$ in $S$ for $\nu\in D^1$, then
$$
F_{pp}^1(0,0)\leq0,\;\;\;F_{pp}^1(0,0)F_{qq}^1(0,0)-(F_{pq}^1)^2(0,0)
\geq0.
$$
Conversely, if
$$
F_{pp}^1(0,0)<0,\;\;\;F_{pp}^1(0,0)F_{qq}^1(0,0)-(F_{pq}^1)^2(0,0)>0,
$$
then $K$ locally maximizes $\R L^1(\nu;f)$ in $S$.
Direct calculation gives
$$\align
H_{uu}(t,x^0,\Psi^1,\pi)&=-2(16t^2-8t+\nu+2), \tag56\\
H_{ux_2}(t,x^0,\Psi^1,\pi)&=4(t+1),\;\;\;H_{ux_4}(t,x^0,\Psi^1,\pi)=4,\tag57\\
H_{u\Psi_2}(t,x^0,\Psi^1,\pi)&=-2,\;\;\;H_{u\Psi_4}(t,x^0,\Psi^1,\pi)=
4(1-3t).\tag58\endalign
$$
The formula (26) remains true as well as $(\Psi_4)_p=2(x_2)_p$ and
$(\Psi_4)_q=2(x_2)_q+1$. Substituting the last relations and (56--58) into
(19) and (20) we obtain
$$\align
u_p&=\frac{(3-5t)2y_4+2y_5-y_6}{16t^2-8t+\nu+2},\tag59\\
u_q&=\frac{(3-5t)2y_{10}+2y_{11}-y_{12}+2(1-3t)}{16t^2-8t+\nu+2}.\tag60\endalign
$$
Evidently, the formulae (29--33), (35--39), and (41--43) are valid for our case,
whereas the equations (34) and (40) are transformed  into
$$\align
\dot y_6&=-4(t+1)u_p-4y_4,\;\;\;y_6(0)=1,\tag61\\
\dot y_{12}&=-4(t+1)u_q-4y_{10}-4,\;\;\;y_{12}(0)=0,\tag62\endalign
$$
respectively.

Summing up above calculation for the evaluation algorithm we state that
by  analogy with Theorems 6 and 7 the problem of finding $\sigma_{24}$ is
reduced to the solution of the equations
$$
\nu y_1(1)+y_3(1)=0\tag63
$$
or
$$
(\nu y_1(1)+y_3(1))(\nu y_7(1)+y_9(1))-(\nu y_{13}(1)+y_{15}(1))^2=0\tag64
$$
where $y_1(t),\dots,y_{15}(t)$ are solutions to the Cauchy problem
for the differential equations (29--33), (61), (35--39), (62), and
(41--43) with $u_p$ and $u_q$ given by (59) and (60). We formulate
the following result.

\proclaim{Theorem 10} Bombieri's number $\sigma_{24}$ is equal to
 the maximum of negative roots of the equations (63) and
(64) multiplied by $(-1)$, where $y_1(t),\dots,y_{15}(t)$ are
solutions to (29--33), (61), (35--39), (62), (41--43), and
(59--60).
\endproclaim

Similarly to Bombieri's number $\sigma_{42}$, the numerical
approximation for $\sigma_{24}$ is $0.969556...$, which is the
maximal negative root of (64).
\medskip

In order to evaluate $\sigma_{34}$ we  consider the functional
$N(\mu;f)=a_4+\mu a_3$ as soon as $\sigma_{34}$ is equal to the supremum
over all real values $\lambda_{34}$ for which $\R N(-\lambda_{34};f)$ is
locally maximized by the Koebe function $K$ in $S$. Now the systems (4) and
(7) are supplied with the transversality conditions
$$
\Psi_3(1)=\mu,\;\;\;\Psi_5(1)=1,\;\;\;\Psi_1(1)=\Psi_2(1)=\Psi_4(1)=0.\tag65
$$
For $u=\pi$ we denote the integrals of (4) and (7)  by $(x^0(t),\Psi^2(t))$,
where $x_1^0(t),\dots,x_4^0(t)$ are given by (12) and
$$\aligned
\Psi_1^2(t)&=t^2-(10+4\mu)t+9+4\mu,\;\;\;\Psi_3^2(t)=-4(t-1)+\mu,\\
\Psi_5^2(t)&=1,\;\;\;\Psi_2^2(t)=\Psi_4^2(t)=0.
\endaligned\tag 66$$
We consider only $\mu\in D^2$, where $D^2$ is a maximal interval on the
$\mu$-axis which satisfies the conditions:
\roster
\item"[i]" $H^2(t,u)=H(t,x^0(t),\Psi^2(t),u)$, as a function of $y=\cos u$, attains
its maximum in $[-1,1]$ only at $y=-1$ for all $t\in[0,1]$;
\item"[ii]" $H_{uu}^2(t,\pi)\neq0,\;\;\;0\leq t\leq1.$
\endroster
Note that $D^2=\{\mu:\mu>-2(\sqrt2-1)\}$.

Again we solve the systems (4) and (7) with $x(0)=0$ and
$$
\Psi_1(0)=9+4\mu,\;\;\;\Psi_2(0)=p,\;\;\;\Psi_3(0)=4+\mu,\;\;\;\Psi_4(0)=q,
\;\;\;\Psi_5(0)=1,\tag67
$$
where $p$ and $q$ are arbitrary real numbers. The conditions (67) are thought of as
variations of (66) at $t=0$.

Let $$ F^2:(p,q)\to\mu x_2(1)+x_3(1) $$ be a real valued mapping
from $(p,q)$-plane into the linear combination of components of
the solution to the Cauchy problem for (4) and (7) with the
initial conditions (67) and $u$ satisfying the maximum principle.
We preserve preceding denotations for $u(t,p,q)$, and similarly to
Theorem 5, for $\mu\in D^2$, we assert, that if $K$ locally
maximizes $\R N(\mu;f)$ in the class $S$, then $$
F_{pp}^2(0,0)\leq0,\;\;\;F_{pp}^2(0,0)F_{qq}^2(0,0)-(F_{pq}^2)^2(0,0)\geq0.
$$ Conversely, if $$
F_{pp}^2(0,0)<0,\;\;\;F_{pp}^2(0,0)F_{qq}^2(0,0)-(F_{pq}^2)^2(0,0)>0,
$$ then $K$ locally maximizes $\R N(\mu;f)$ in $S$.

Direct calculation gives
$$\align
H_{uu}(t,x^0,\Psi^2,\pi)&=-4(8t^2-(4-2\mu)t+1),\tag68\\
H_{ux_2}(t,x^0,\Psi^2,\pi)&=4(t+1+\mu),\;\;\;H_{ux_4}(t,x^0,\Psi^2,\pi)=4,
\tag69\\
H_{u\Psi_2}(t,x^0,\Psi^2,\pi)&=-2,\;\;\;H_{u\Psi_4}(t,x^0,\Psi^2,\pi)=
4(1-3t).\tag70\endalign
$$
The formulae (19), (20), and (68--70) lead to
$$\align
u_p&=\frac{(3-5t+\mu)2y_4+2y_5-y_6}{16t^2-(8-4\mu)t+2}\tag71\\
u_q&=\frac{(3-5t+\mu)2y_{10}+2y_{11}-y_{12}+2(1-3t)}{16t^2-(8-4\mu)t+2}.
\tag72\endalign
$$
The formulae (34) and (40) are changed to
$$\align
\dot y_6&=-4(t+1+\mu)u_p-4y_4,\;\;\;y_6(0)=1,\tag73\\
\dot y_{12}&=-4(t+1+\mu)u_q-4y_{10}-4,\;\;\;y_{12}(0)=0,\tag74\endalign
$$
respectively.

Let us sum up the results and state, that similarly to Theorems 6 and 7, the
problem of finding $\sigma_{34}$ is reduced to the solution of the equations
$$
\mu y_2(1)+y_3(1)=0,\tag75
$$
or
$$
(\mu y_2(1)+y_3(1))(\mu y_8(1)+y_9(1))-(\mu y_{14}(1)+y_{15}(1))^2=0,\tag76
$$
where $y_1(t),\dots,y_{15}(t)$ are solutions to the Cauchy problem for
the differential equations (29--33), (73), (35--39), (74), and (41--43)
with $u_p$ and $u_q$ given by (71) and (72).

\proclaim{Theorem 11} Bombieri's number $\sigma_{34}$ is equal to
the maximum of the negative roots of equations (75) and (76)
multiplied by $(-1)$, where $y_1(t),\dots,y_{15}(t)$ are solutions
of (29--33), (73), (35--39), (74), (41--43), and
(71--72).\endproclaim

Numerical methods applied to (75) and (76) show that
$\sigma_{34}\approx 0.791557...$

\head \S 7. Fourth coefficient of bounded univalent functions\endhead

In this section we will find a sharp constant $M_0$, such that for
all $M\geq M_0$ the Pick function $P_M$ gives the local maximum to
$\R a_4$ (and to $|a_4|$) over all univalent functions $f\in S^M$.
In the class $S_R^M\subset S^M$ of functions with real Taylor
coefficients the Pick function gives the maximum to $a_4$ for
$M\geq 11$ and this constant is sharp (see [17], [19, p.163]).
That is why we can consider only $M\geq 11$. For this problem an
analog of Theorem 5 holds and allows us to construct an evaluation
algorithm. Therefore, we again consider the systems (4) and (7)
supplied now with the transversality conditions $$
\Psi_5(1-1/M)=1,\;\;\;\Psi_1(1-1/M)=\dots=\Psi_4(1-1/M)=0.\tag77
$$ For $u=\pi$, let us denote the integrals of (4) and (7)  by
$(x^0(t),\Psi^3(t))$, where $$\aligned
\Psi_1^3(t)&=t^2-\left(10-\frac{16}{M}\right)t+9-\frac{24}{M}+\frac{15}{M^2}
, \;\;\;\Psi_3^3(t)=-4\left(t-1+\frac{1}{M}\right),\\
\Psi_5^3(t)&=1,\;\;\;\Psi_2^3(t)=\Psi_4^3(t)=0.
\endaligned\tag78$$
Varying conditions (78) at $t=0$ we put
$$
\Psi_1(0)=9-\frac{24}{M}+\frac{15}{M^2},\;\;\;\Psi_2(0)=p,\;\;\;
\Psi_3(0)=4\left(1-\frac{1}{M}\right),\;\;\;\Psi_4(0)=q,\;\;\;\Psi_5(0)=1,
\tag79
$$
where $(p,q)\in\Bbb R^2$.

Let
$$
F^3:(p,q)\to x_5(1-1/M)
$$
be a mapping from the $(p,q)$-plane into the axis of the third component of the solution to the
Cauchy problem for (4) and (7) with the initial conditions (79) and $u$
satisfying the maximum principle. Preserving notations of Section 6 for
$u(t,p,q)$, similarly to Theorem 5, we observe that for $M>11$, if
$P_M$ locally maximizes $\R a_4$ in $S(M)$, then
$$
F_{pp}^3(0,0)\leq0,\;\;\;F_{pp}^3(0,0)F_{qq}^3(0,0)-(F_{pq}^3)^2(0,0)
\geq0,
$$
and conversely, if
$$
F_{pp}^3(0,0)<0,\;\;\;F_{pp}^3(0,0)F_{qq}^3(0,0)-(F_{pq}^3)^2(0,0)>0,
$$
then $P_M$ locally maximizes $\R a_4$ in $S(M)$.

Direct calculation gives
$$\align
H_{uu}(t,x^0,\Psi^3,\pi)&=-2\left(16t^2-\left(8+\frac{16}{M}\right)t+2-
\frac{8}{M}+\frac{15}{M^2}\right),\tag80\\
H_{ux_2}(t,x^0,\Psi^3,\pi)&=4\left(t+1-\frac{4}{M}\right),\;\;\;
H_{ux_4}(t,x^0,\Psi^3.\pi)=4,\tag81\\
H_{u\Psi_2}(t,x^0,\Psi^3,\pi)&=-2,\;\;\;
H_{u\Psi_4}(t,x^0,\Psi^3,\pi)=4(1-3t).\tag82\endalign
$$
The formulae (19), (20), and (80--82) lead to
$$\align
u_p&=\frac{(3-4/M-5t)2y_4+2y_5-y_6}{16t^2-(8+16/M)t+2-8/M+15/M^2},\tag83\\
u_q&=\frac{(3-4/M-5t)2y_{10}+2y_{11}-y_{12}+2(1-3t)}{16t^2-(8+16/M)t+2-
8/M+15/M^2}.\tag84\endalign
$$
The formulae (34) and (40) are changed to
$$\align
\dot y_6&=-4\left(t+1-\frac{4}{M}\right)u_p-4y_4,\;\;\;y_6(0)=1,\tag85\\
\dot y_{12}&=-4\left(t+1-\frac{4}{M}\right)u_q-4y_{10}-4,\;\;\;y_{12}(0)=0,
\tag86
\endalign$$
respectively.

Summing up the results we state, that similarly to Theorems 6 and 7, the
problem of finding the best possible $M$ for which the Pick function $P_M$
locally maximizes $\R a_4$ in $S(M)$ is reduced to the solution of the equations
$$
y_3(1-1/M)=0,\tag87
$$
or
$$
y_3(1-1/M)y_9(1-1/M)-y_{15}^2(1-1/M)=0,\tag88
$$
where $y_1(t),\dots,y_{15}(t)$ are solutions to the Cauchy problem for
the differential equations (29--33), (85), (35--39), (86), and (41--43)
with $u_p$ and $u_q$ given by (83) and (84).

\proclaim{Theorem 12} The Pick function $P_M$ locally maximizes
$\R a_4$ in $S(M)$ if $M>M_0$ and does not give a local maximum to
$\R a_4$ if $M<M_0$, where $M_0$ is the maximum of the roots of
the equations (87) and (88), where $y_1(t),\dots,y_{15}(t)$ are
solutions of (29--33), (85), (35--39), (86), (41--43), and
(83--84).\endproclaim

Numerical methods applied to (87) and (88) show that $M_0\approx
22.9569...$

\enddocument